\documentclass[12pt]{article}

\usepackage[english]{babel}
\usepackage{amssymb}
\usepackage{graphics}

\newcommand{\qed}{$\;\;\;\Box$}

\newcounter{claim}[section]
\renewcommand{\theclaim}{\arabic{claim}}
{\par\medskip\par}

\newcommand{\A}{{\cal A}}                       
\newcommand{\N}{{\mathbb N}}
\newcommand{\R}{{\mathbb R}}
                       
\newcommand{\diag}{\mathop{\rm diag}\nolimits}

\newcommand{\hide}[1]{}
\newtheorem{thm}{Theorem}[section]

\newtheorem{defn}[thm]{Definition}

{\qed\par\smallbreak}

\date{}
\title{1D Hyperbolic Systems with  \\ Nonlinear Boundary Conditions II: \\ Criteria for Finite Time Stability}

\newcounter{thesame}
\author{
	Irina Kmit
	\thanks{Institute of Mathematics, Humboldt University of Berlin. On leave from the
		Institute for Applied Problems of Mechanics and Mathematics, 
		National Academy of Sciences of	Ukraine, Lviv, Ukraine. 
		{\small   E-mail:
			{\tt kmit@mathematik.hu-berlin.de}}}
}

\begin{document}
	\maketitle

	\begin{abstract}
		\noindent
		We investigate the finite time stability property of one-dimensional nonautonomous initial  boundary value  problems for  linear  decoupled hyperbolic systems with nonlinear  boundary conditions.
		We  establish sufficient and necessary conditions 	
		under which  continuous or $L^2$-generalized solutions stabilize to zero in a finite time. 	Our criteria are expressed in terms 	of a propagation operator along characteristic curves.
	
	\end{abstract}

\section{Introduction}
\renewcommand{\theequation}{{\thesection}.\arabic{equation}}
\setcounter{equation}{0}
\subsection{Problem}\label{sec:problem}

Established in the middle of the 50th, the 
Finite Time Stability  (FTS)
concept attracts  growing attention in view of its applications  
in  control and system engineering
\cite{corbas,gugat,lakra,per,udw05,udw12},
output-feedback stabilization 
\cite{hege,hegez,BerKrs,xu}, inverse problems
 \cite{Pril, tikh18}), ATM networks \cite{AbdCos02}, car suspension systems \cite{CosTom14}, and robot manipulators \cite{Or105}.  This concept is used in two ways.
Quantitatively, it 
describes a restrained behavior of the dynamical system over a specified time interval. Qualitatively, it
characterizes asymptotically stable dynamical systems whose
trajectories reach an equilibrium point in a finite time.
In this paper we characterize FTS hyperbolic systems using the qualitative notion of FTS. 

In \cite{KL2} we gave a
comprehensive FTS analysis of a class of linear initial-boundary value
 problems with reflection boundary conditions for  decoupled nonautonomous hyperbolic systems, providing  algebraic and combinatorial
 criteria.
 In the autonomous setting, we provided also a spectral criterion. Asymptotic properties of solutions to perturbed FTS problems were studied in \cite{KmLyul}.  In the present paper, we  establish  FTS criteria for a class of nonlinear boundary value problems. These results can be
applied to solving inverse problems for  hyperbolic systems with FTS boundary conditions (as we demonstrate in Subsection~\ref{inverse}).

Let $n\ge 2$. 
Our stability results   concern the decoupled nonautonomous hyperbolic system
\begin{equation}\label{ss2d}
	\partial_tu + A(x,t)\partial_xu + B(x,t)u = 0, \quad 0<x<1, t>0,
\end{equation}
where $u=(u_1,\ldots,u_n)$   is a vector of real-valued functions and the diagonal matrices $A=\diag(a_1,\dots,a_n)$ 
and $B=\diag(b_{1},\dots,b_{n})$  have real entries. 

 Set 
$
\Pi=\{(x,t):\,0\le x\le 1,\, t\ge 0\}.
$ 
Suppose that 
\begin{equation}\label{ss4}
	\inf_{(x,t)\in \Pi}a_j\ge a \ \mbox{ for all }\  j\le m\quad 
	\mbox{ and }\quad 
	\sup_{(x,t)\in \Pi}a_j\le -a \ \mbox{ for all } \ j>m
\end{equation}
for some $a>0$ and $0\le m\le n$. The system (\ref{ss2d})  is 
subjected to the initial conditions
\begin{equation}\label{ss11}
	u(x,0)=\varphi(x), \quad 0\le x\le 1,
\end{equation}
and the homogeneous nonlinear boundary conditions 
\begin{equation}\label{ss5}
	u^{out}(t)=h(t,u^{in}(t)), \quad t\ge 0,
\end{equation}
where $h=h(t,\xi)=(h_1(t,\xi),\dots,h_n(t,\xi))$, with   $\xi\in \R^n$,
is a real valued function,
\begin{equation}\label{ss51}
	h(t,0)=0\quad \mbox{for all } t\ge 0,
\end{equation}
and
$$
	\begin{array}{rcl}
		u^{out}(t)&=&(u_1(0,t),\dots,u_m(0,t),u_{m+1}(1,t),\dots,u_n(1,t)),\\ [1mm]
		u^{in}(t)&=&(u_1(1,t),\dots,u_m(1,t),u_{m+1}(0,t),\dots,u_n(0,t)).
	\end{array}
$$

\subsection{Preliminaries on continuous and $L^2$-generalized solutions}

Let
\begin{equation}\label{v1}
	\begin{array}{rcl}
		\varphi^{out}&=& (\varphi_1(0),\dots,\varphi_m(0),\varphi_{m+1}(1),\dots,\varphi_n(1)),\\ [1mm] \varphi^{in}&=&(\varphi_1(1),\dots,\varphi_m(1),\varphi_{m+1}(0),\dots,\varphi_n(0)). 
	\end{array}
\end{equation}
We say that a function  $\varphi$ satisfies the zero order compatibility conditions between  (\ref{ss11}) and (\ref{ss5}) if
\begin{equation}\label{ss50}
	\varphi^{out}=h(0,\varphi^{in}).
\end{equation}
We consider the 
set $C_h(\Pi)^n$ of functions   $u\in C(\Pi)^n$ such that 
$
u^{out}(0)=h(0,u^{in}(0)).
$
Note that, if $u\in C_h(\Pi)^n$, then  $u(x,0)$ satisfies the zero order
compatibility conditions between (\ref{ss11}) and (\ref{ss5}) with $\varphi=u(x,0)$. 
Let $C_h([0,1])^n$ be
a closed subset of a Banach space  $C([0,1])^n$ that  consists of functions
$\varphi\in C([0,1])^n$ fulfilling the condition (\ref{ss50}). 
Furthermore,  $C_h^1([0,1])^n= C_h([0,1])^n\cap C^1([0,1])^n$.

Let us introduce  solution concepts, that will be  used in the paper.
To this end, we first define characteristics of (\ref{ss2d}) as follows.
 For given $j\le n$,
$x \in [0,1]$, and $t >0$, 
the $j$-th characteristic of (\ref{ss2d})
passing through the point $(x,t)\in\Pi$ is
the solution $\omega_j(\xi)=\omega_j(\xi,x,t): [0,1]\to\R$ to the initial value problem 
$$
\partial_\xi\omega_j(\xi, x,t)=\frac{1}{a_j(\xi,\omega_j(\xi,x,t))},\;\;
\omega_j(x,x,t)=t.
$$

Let a continuous function $u:\Pi\to \R^n$ be  continuously differentiable  in~$\Pi$ excepting at most a countable number of
characteristic curves of  (\ref{ss2d}). If $u$  satisfies
(\ref{ss2d}), (\ref{ss11}), and (\ref{ss5})  in $\Pi$ except the aforementioned characteristic curves, then it
   is called a  {\it piecewise continuously differentiable  solution}
 to the problem (\ref{ss2d}), (\ref{ss11}),~(\ref{ss5}).

 If the initial function $\varphi$ is sufficiently smooth, then  using integration along characteristics, we can transform  the problem (\ref{ss2d}), (\ref{ss11}), (\ref{ss5})
to a system of  integral equations.
The characteristic curve $\tau=\omega_j(\xi,x,t)$ reaches the
boundary of $\Pi$ in two points with distinct ordinates. Let $x_j(x,t)$
denote the abscissa of that point whose ordinate is smaller.
Note that  the value of  $x_j(x,t)$
does not depend on $x$ and $t$ if $t>1/a$, where
$a>0$ satisfies~(\ref{ss4}). More precisely, if $t>1/a$, then
\begin{eqnarray*}\label{xj}
	x_j(x,t)=x_j=\left\{
	\begin{array}{rl}
		0 &\mbox{if}\ 1\le j\le m\\
		1 &\mbox{if}\ m<j\le n.
	\end{array}
	\right.
\end{eqnarray*}
Set
$$
c_j(\xi,x,t)=\exp \int_x^\xi
\left(\frac{b_{j}}{a_{j}}\right)(\eta,\omega_j(\eta,x,t))\,d\eta.
$$
Define  a linear operator
$S: C(R_{+})^n\to C(\Pi)^n$
by 
$$
[Sv]_j(x,t)=c_j(x_j(x,t),x,t)v_j(\omega_j(x_j(x,t),x,t)), \quad j\le n,
$$
and a nonlinear operator 
$R:C(\Pi)^n\to C(R_{+})^n$ by
\begin{eqnarray*}
	\displaystyle
	\left[Ru\right]_j(t) = h_j(t,u^{in}(t)),
	\quad   j\le n.
\end{eqnarray*}

As it follows from the method of characteristics,
any piecewise continuously differentiable  solution $u$
 to the problem (\ref{ss2d}), (\ref{ss11}), (\ref{ss5}) 
 satisfies  the following system of functional equations:
\begin{equation}
	\begin{array}{rr}
\label{rep1}
u_j(x,t)=[Qu]_j(x,t)
	\end{array}
\end{equation}
where the affine operator $Q: D(Q)\subset  C_h(\Pi)^n\to C_h(\Pi)^n$  is defined by
\begin{equation}
\label{Q}
[Qu]_j(x,t)=
\left\{\begin{array}{lcl}
\displaystyle \left[SRu\right]_j(x,t)
& \mbox{ if }&  x_j(x,t)=0 \mbox{ or } x_j(x,t)=1 \\
c_j(x_j(x,t),x,t)\varphi_j(x_j(x,t))      & \mbox{ if }&  x_j(x,t)\in(0,1),
\end{array}\right.
\end{equation}
and
$$
D(Q)=\{u\in C_h(\Pi)^n\,:\, u(x,0)=\varphi(x)\}.
$$
Note that the definition of $Q$ depends on the choice of the
function $\varphi$. We will write $Q=Q_\varphi$ when we want to
specify this dependence explicitly.

Vice versa, if a $C$-map $u:{\Pi}\to \R^n$   is
piecewise continuously differentiable 
excepting at most a countable number of
characteristic curves of  (\ref{ss2d})
and satisfies~(\ref{rep1}) pointwise, then it
is a piecewise continuously differentiable solution to 
(\ref{ss2d}), (\ref{ss11}), (\ref{ss5}). This motivates the following definition.

\begin{defn}
A continuous function~$u:\Pi\to\R^n$ satisfying  (\ref{rep1}) in $\Pi$
is called a {\it continuous solution} to (\ref{ss2d}), (\ref{ss11}), (\ref{ss5}).
\end{defn}

For a Banach space $X$, the $n$-th Cartesian power $X^n$ is considered to be a Banach space of vectors $u=(u_1,\dots,u_n)$ normed by   $
\|u\|_{X^n}= \max_{i\le n} \|u_i\|_{X}.
$ 
Let $\|\cdot\|_{\max}=\max_{jk}|m_{jk}|$  denote the $\max$-matrix norm of $M=(m_{jk})$   in the space of matrices
$M_n$.

Below  we will use  our result from  \cite[Theorem 3.1]{ijdsde} about the existence and uniqueness
of global regular solutions. 
 \begin{thm}\label{thm:classical} 
 Let  the condition
 (\ref{ss4}) 
  be
 fulfilled. Moreover, 
 assume that 
 	\begin{equation}\label{ss53}
 	\begin{array}{ll}
 		\mbox{for all }  j,k\le n \mbox{ the functions } a_j,   b_{j}, 
 		\mbox{ and } h_j \\\mbox{are continuously differentiable in all their arguments}
 \end{array}\end{equation}
and  
 for each $T>0$
 there exists a positive real $C(T)$ and a 	polynomial $H$ such that
	\begin{equation}\label{eq:h}
	\left\{\left\|\nabla_\xi h(t,\xi)\right\|_{max}\,:\, 0\le t\le T,\ \xi\in \R^n\right\}\le C(T)\left(\log\log H(\|\xi\|)\right)^{1/4}.
\end{equation}
Then the following is true.

{\bf 1.}  For every $\varphi\in C_h([0,1])^n$,
the problem~(\ref{ss2d}), (\ref{ss11}), (\ref{ss5}) has a unique   continuous solution in
$\Pi$.

{\bf 2.}
For every $\varphi\in C^1_h([0,1])^n$,
the problem~(\ref{ss2d}), (\ref{ss11}), (\ref{ss5}) has a unique piecewise continuously differentiable solution in
$\Pi$.
\end{thm}

We now define an $L^2$-generalized solution 
to the problem  (\ref{ss2d}), (\ref{ss11}), (\ref{ss5}) similarly to
\cite[Definition 2]{lyul_I}.

\begin{defn} \rm\label{defn:sol}
Assume that the conditions of Theorem 
	\ref{thm:classical} are fulfilled. 	Let  $\varphi\in L^2(0,1)^n$.
	A function $u\in C\left([0,\infty), L^2(0,1)\right)^n$ is called
	an {\it $L^2$-generalized  solution} to the problem (\ref{ss2d}), (\ref{ss11}), (\ref{ss5}) 
	if, for any sequence $\varphi^l\in C^1_{h}([0,1])^n$ with
	$\varphi^l$ converging to $\varphi$ in $L^2(0,1)^n$,
	the sequence  of piecewise continuously differentiable solutions
	$u^l(x,t)$ to the problem
	(\ref{ss2d}), (\ref{ss11}), (\ref{ss5}) with
	$\varphi$ replaced by $\varphi^l$  fulfills the convergence condition
	\begin{equation}\label{con}
	\|u^l(\cdot,t)-u(\cdot,t)\|_{L^2(0,1)^n} \to 0\quad\mbox{as}\quad
	l\to\infty,
	\end{equation}
	uniformly in $t$ varying in the range $0\le t\le T$, for each $T>0$.
\end{defn}
Here the norm in $L^2(0,1)^n$ is defined as usual by
$
\|u\|_{L^2(0,1)^n}^2=\int_0^1(u,u)\,dx= \int_0^1 \sum_{i=1}^n u_i^2\, dx,
$
where  $(\cdot,\cdot)$ here and below denotes the scalar product  in $\R^n$.

The following existence and uniqueness result is obtained in \cite[Theorem 2]{lyul_I}.

\begin{thm}\label{evol} 
		Let  the conditions 	(\ref{ss4}), 	(\ref{ss51}), and (\ref{ss53}) be 	fulfilled. 
Moreover, assume that  for each $T>0$
	there exists a positive real $C(T)$ such that
	\begin{equation}\label{eq:h1}
		\sup\left\{\|\nabla_\xi h(t,\xi)\|_{\max}\,:\, 0 \le t\le T,\ \xi\in \R^n\right\}\le C(T).
	\end{equation} 
	Then, for every $\varphi\in L^2(0,1)^n$,  the problem
	(\ref{ss2d}), (\ref{ss11}), (\ref{ss5}) has a unique
	$L^2$-generalized  solution.
\end{thm}

\subsection{Our results}\label{sec:problem}

If the  problem (\ref{ss2d}), (\ref{ss11}), (\ref{ss5}),  (\ref{eq:h})
has  an $L^2$-generalized solution, then it is unique just by
Definition \ref{defn:sol}. If
this  problem 
has  a continuous solution, it is also unique as shown in 
 \cite{ijdsde}
(see the proof of
 \cite[Theorem 3.1]{ijdsde}).

\begin{defn}\label{defn:stab}\rm
Assume that,	for every $\varphi\in L^2(0,1)^n$ (resp.,  $\varphi\in C_h([0,1])^n$), the problem (\ref{ss2d}), (\ref{ss11}), (\ref{ss5}),  (\ref{eq:h})  has an $L^2$-generalized  solution (resp., a continuous
	solution). We say that this
	 problem   is    {\it Finite Time Stabilizable (FTS)} 
if
there exists  a positive real $T$  such that,
for every $\varphi\in L^2(0,1)^n$ (resp.,  $\varphi\in C_h([0,1])^n$), the $L^2$-generalized  solution (resp., a continuous
solution)  
   is a constant zero function for   $t>T$. 
  The infimum  of all  $T$ with the
above property is called the {\it optimal  stabilization time} and is denoted by  $T_{opt}$.
\end{defn}

Since  the operator
 $Q$ operates with functions on shifted domains
 and, thus, captures propagation  from the boundary $\partial\Pi$ into the domain $\Pi$,  the stabilization properties heavily depend on the powers of the operator $Q$. We start with a useful property of the  operator~$Q$. Given
  $T>0$, set $
\Pi^T=\{(x,t)\in\Pi\,:\,t\le T\}.$
\begin{thm}\label{control_inverse}
	For  every $T>0$ there exists $k\in \N$ such that the following is true. If, for $w\in C_h(\Pi)^n$, the 
	problem 
	(\ref{ss2d}), (\ref{ss11}), (\ref{ss5}),  (\ref{eq:h}) 
	with $\varphi(x)=w(x,0)$
	has a unique continuous solution $u$ in $\Pi$, then 
	 $u(x,t)=[Q^kw](x,t)$ in $\Pi^T$ where 
	$Q=Q_\varphi$ for  $\varphi(x)=w(x,0)$.
	  \end{thm}

Now we formulate  our stabilization criterion in the nonautonomous setting.
 \begin{thm}\label{tm_1}
Let the condition (\ref{ss51})
 be fulfilled.  Assume that,	for every $\varphi\in L^2(0,1)^n$ (resp.,  $\varphi\in C_h([0,1])^n$), the problem (\ref{ss2d}), (\ref{ss11}), (\ref{ss5}),  (\ref{eq:h})  has an $L^2$-generalized  solution (resp., a continuous
 solution).  
 Then this problem
 is FTS
if and only if
\begin{equation}
\label{C0}
\begin{array}{ll}
\mbox{ there is }
T>0 \mbox{ and }
k\in \N    \mbox{ such that,  }
\mbox{ for all }
w\in C_h(\Pi)^n \mbox{ and } x\in[0,1],\\
\left[Q^kw\right]\left(x,T\right)\equiv 0 \mbox{ where }
Q=Q_\varphi \mbox{ for } \varphi(x)=w(x,0).
\end{array}
\end{equation}
\end{thm}

In the autonomous setting a stabilization  criterion is formulated  
 in a stronger form. 

\begin{thm}\label{tm_1_autonom} 
		Assume that  the coefficient matrices $A$ and $B$ do not depend on $t$ and the boundary function $h$ does not explicitely depend on $t$, that is, $h(t,\xi)\equiv h(\xi)$.
	Moreover, let the condition  (\ref{ss51}) 
		be fulfilled. Assume also that,	for every $\varphi\in L^2(0,1)^n$ (resp., $\varphi\in C_h([0,1])^n$), the problem (\ref{ss2d}), (\ref{ss11}), (\ref{ss5}),  (\ref{eq:h})  has an $L^2$-generalized  solution (resp., a continuous
		solution). 		
			Then this problem
  is FTS 
		if and only if 
	\begin{equation}\label{C00}
	\begin{array}{ll}
		\mbox{ there is }
		T>0 \mbox{ and }
		q\in \N    \mbox{ such that,   for all } k\in \N,\
			w\in C_h(\Pi)^n,  \mbox{ and }
		 x\in[0,1],\\
		\left[Q^{kq}w\right]\left(x,kT\right)= 0   
	 \mbox{ where }
	Q=Q_\varphi \mbox{ for } \varphi(x)=w(x,0).
	\end{array}
	\end{equation}
\end{thm}

Theorems \ref{control_inverse}--\ref{tm_1_autonom} assume the existence of $L^2$-generalized or continuous
solutions (recall that those are always unique). While some sufficient conditions for the existence of
solutions to the problem (\ref{ss2d}), (\ref{ss11}), (\ref{ss5}), (\ref{eq:h})
are given in Theorems~\ref{thm:classical} and~\ref{evol}, we want to emphasize that
Theorems \ref{control_inverse}--\ref{tm_1_autonom} are not restricted to these particular
conditions and are more general.

The rest of the paper is organized as follows.   
The FTS-criteria  of Theorems    \ref{tm_1} and \ref{tm_1_autonom}  are proved in 
Section~\ref{sec:criteria}.
Discussion of our stabilization criteria  are provided in Section \ref{ex}, where we also show how our 
Theorem~\ref{control_inverse} can be applied to solving
inverse hyperbolic problems.

\section{Stabilization criteria}\label{sec:criteria}
\renewcommand{\theequation}{{\thesection}.\arabic{equation}}
\setcounter{equation}{0}
\subsection{ Proof of 
	Theorem~\ref{control_inverse}}

Fix an arbitrary $T>0$. 
Since $Q$ is a down-shift operator along characteristic curves up to the boundary of $\Pi$
in the direction of time decrease,
there exists an integer $q=q(T)$
such that all iterations of the operator $Q$ starting from the $q$-th iteration 
stabilize, namely for every $w\in C_h(\Pi)^n$ it holds in $\Pi^T$ that
\begin{equation}\label{tQ0}
\left[Q^qw\right](x,t)=\left[ Q^{q+1}w\right](x,t),
\end{equation} 
where in the definition  (\ref{Q}) of the operator $Q$ we set 
 $\varphi(x)=w(x,0)$.

Fix a function $w\in C_h(\Pi)^n$ fulfilling the conditions of
Theorem~\ref{control_inverse}.
Then the problem (\ref{ss2d}), (\ref{ss11}), (\ref{ss5}), (\ref{eq:h}) with $\varphi=w(x,0)$ has a unique continuous solution.
 Set $u=Q^qw$. Hence,   $u\in C_h(\Pi)^n$, and (\ref{tQ0}) implies that in $\Pi^T$ we have
$$
[Qu](x,t)=[Q^{q+1}w](x,t)=[Q^{q}w](x,t)=u(x,t).
$$
It follows that the function $u=Q^{q}w$ is the continuous solution in $\Pi^T$ to the problem 
(\ref{ss2d}), (\ref{ss11}), (\ref{ss5}), (\ref{eq:h}) with $\varphi=w(x,0)$. 
The proof of Theorem \ref{control_inverse} is complete. 

\subsection{Nonautonomous case: proof of Theorem \ref{tm_1}} \label{sec:prcrit}

{\it Sufficiency.} Let $T>0$ and $k\in\N$ be numbers satisfying  the 
condition~(\ref{C0}).
Fix an arbitrary $\varphi\in L^2(0,1)^n$. Suppose that
 the problem (\ref{ss2d}), (\ref{ss11}), (\ref{ss5}),  (\ref{eq:h})
 has a unique
$L^2$-generalized  solution $u$. 

First note that $C^1_{h}([0,1])^n $ is densely embedded  into $L^2(0,1)^n.$
Indeed, since the boundary conditions  (\ref{ss5}) are homogeneous (see \ref{ss51}),
$C_0^\infty([0,1])^n$   is a subset of 
$C^1_{h}([0,1])^n$.  As ususal, by $C_0^\infty([0,1])$ we denote a subspace of  $C^\infty([0,1])$
that consists  of
functions having  support within $(0,1)$. 
Now, we  fix an arbitrary sequence $\varphi^l\in C_h^1([0,1])^n$ such that
$\varphi^l$ converges to $\varphi$ in $L^2(0,1)^n$ 
 and let $u^l(x,t)$ be the piecewise continuously
differentiable
 solution to the problem
(\ref{ss2d}), (\ref{ss11}), (\ref{ss5}), (\ref{eq:h}) with  $\varphi$ replaced by~$\varphi^l$ (see Theorem \ref{thm:classical}).

By Definition  \ref{defn:sol}, the sequence $u^l(x,t)$
converges as in~(\ref{con}).
Using integration along characteristics,
we see that
$$
u^l(x,t)=[Qu^l](x,t)
\quad \mbox{ for all } x\in[0,1] \mbox{ and } t\in[0,T].
$$
This means that  the function $u^l(x,t)$ is a fixed point of the operator
$Q$
and, hence, of any power of $Q$.
Combining this with the condition (\ref{C0}), 
we conclude that
$$
u^l(x,T)=
\left[Q^ku^l\right](x,T)=0 \quad \mbox{for all } x\in[0,1] \mbox{ and } l\in \N.
$$
Since the initial boundary value problem (\ref{ss2d}),  (\ref{ss5}), (\ref{eq:h}) with the zero initial data  at  $t=T$ has a unique
 piecewise continuously differentiable  solution for
 $t\ge T$ 
(see Theorem \ref{thm:classical}), we conclude
that  $u^l\equiv 0$ for
 $t\ge T$ . The identity  $u\equiv 0$ for $t>T$ follows from the convergence (\ref{con}).
The  FTS property is therewith proved.

If the problem (\ref{ss2d}), (\ref{ss11}), (\ref{ss5}),  (\ref{eq:h})
has a unique continuous solution, the proof goes along the same lines as above with obvious simlifications.

{\it Necessity.}
Consider  first the case when the problem (\ref{ss2d}), (\ref{ss11}), (\ref{ss5}),  (\ref{eq:h}) is FTS and all $L^2$-generalized solutions stabilize to
zero in a finite time.
Fix an arbitrary $T>T_{opt}$ and an integer $q=q(T)$ fulfilling the condition (\ref{tQ0}) in $\Pi^T$. 
Fix an arbitrary $w\in C_h(\Pi)^n$ and put $\varphi(x)=w(x,0)\in C_h([0,1])$. Then, by assumption, the 
problem (\ref{ss2d}), (\ref{ss11}), (\ref{ss5}),  (\ref{eq:h}) 
has a unique $L^2$-generalized solution. Moreover, as 
$\varphi\in C_h([0,1])$, then by Theorem \ref{thm:classical}, this problem has a unique continuous solution. We, therefore, fall into the conditions
of Theorem \ref{control_inverse}. As shown in the proof of Theorem \ref{control_inverse}, the function $
u=Q^qw\in C_h(\Pi)^n$ is a continuous solution in $ \Pi^T$   to
the problem (\ref{ss2d}), (\ref{ss11}), (\ref{ss5}),  (\ref{eq:h}). 
Since any continuous solution is an $L^2$-generalized solution, 
then using the FTS property for the $L^2$-generalized solutions,
we conclude that
$
\left[Q^qw\right](x,T)=0$
for all  $x\in[0,1],
$
as desired. 

If the problem (\ref{ss2d}), (\ref{ss11}), (\ref{ss5}),  (\ref{eq:h}) is FTS and all continuous solutions stabilize to zero in a finite time, the argument is similar and even simpler than in the case we considered.

The proof of  Theorem  \ref{tm_1} is complete.

\subsection{Autonomous case: proof of 
	Theorem~\ref{tm_1_autonom}}

\hskip6mm {\it Sufficiency.}
Since the condition (\ref{C00}) implies (\ref{C0}),
 this part immediately follows from the sufficiency part of Theorem \ref{tm_1}.

{\it Necessity.} 
Consider two cases.

{\bf Case 1:} the problem (\ref{ss2d}), (\ref{ss11}), (\ref{ss5}),  (\ref{eq:h}) is FTS and all continuous solutions stabilize to zero in a finite time.
Fix 
$T>T_{opt}$ and  $q\in \N$ fulfilling both
the condition (\ref{C0}) with $k=q$ and the equality
(\ref{tQ0}) in $\Pi^{2T}$.
For any continuous solution $u$
we have 
\begin{equation}\label{q1}
0=\left[Q^{q}u\right](x,t)=\left[(SR)^{q}u\right](x,t)\quad \mbox{for all } x\in[0,1]
\end{equation} 
and for all $t\ge T$, where the second equality can be proved
as follows.
 We first prove that this  equality is fulfilled for all $t\in[T,2T]$. By the way of contradiction, assume that  this 
 is not true for some  continuous solution~$u$. 
	Then there exist $x\in[0,1]$, $t\in[T,2T]$, and $j\le n$ such that 
	the value $\left[ Q^qu\right]_j(x,t)$
		can be  expressed  in terms of the
	values of $u$ at points lying on the initial axis.
		Straightforward calculations show that
	there exist
	positive  integers $q_1,\dots,q_n$. as well as
	 $C^1$-functions $ F: \R^{q_1+\dots +q_n}\mapsto \R$
	 and  $\widetilde F: \R^{q_1}\times\dots\times \R^{q_n}\mapsto \R$, 	and  pairwise distinct reals $x_{sr}\in [0,1]$ such that
		\begin{equation}\label{QQ}
	\begin{array}{cc}
		\left[Q^qu\right]_j(x,t)=\widetilde{F}(\bar v_{1}^u,\dots,\bar v_{n}^u),
	\end{array}
	\end{equation} 
	where 
	$$
	\widetilde{F}(\bar v_{1}^u,\dots,\bar v_{n}^u)={F}(v_{1}^u,
	v_2^u,\dots,v_{q_1}^u, v_{q_1+1}^u,\dots, v_{q_1+q_2}^u,v_{q_1+q_2+1}^u,\dots, v_{q_1+\dots+q_{n}}^u)
	$$
	and the vector-function $\bar v_{s}^u$ for all $s\le n$ is given by
	\begin{equation}\label{baru}
	\bar v_{s}^u=\left(v_{q_1+q_2+\dots+q_{s-1}+1}^u,\dots, v_{q_1+q_2+\dots+q_{s}}^u\right)
	=\left(u_s(x_{s1},0),\dots,u_s(x_{sq_s},0) \right).
	\end{equation} 
	Since $u$ is a solution, we have $\varphi(x)=u(x,0)$. It follows that
	  $\widetilde{F}$ is  a \underline{} composition of two homogeneous operators, namely the multiplication-shift operator $S$
	and the  nonlinear  boundary operator $R$. This implies that 
	$\widetilde {F}(0,\dots,0)= 0$.
		Note that, due to (\ref{tQ0}) in $\Pi^{2T}$, the representation 	(\ref{QQ}) is unique.
	
Equality (\ref{tQ0}) considered in $\Pi^{2T}$ implies that
$u(x,t)=\left[Q^qu\right](x,t)$. Combined with
 (\ref{QQ}), this gives the equality
	\begin{equation}\label{Qphi1}
	\begin{array}{rcl}
		u_j(x,t)&=&\left[Q^qu\right]_j(x,t)=	\widetilde{F}(\bar v_{1}^u,\dots,\bar v_{n}^u)=\displaystyle
		\widetilde{F}(\bar v_{1}^u,\dots,\bar v_{n}^u)-\widetilde{F}(0,\dots,0)\\
		&=&\displaystyle
		\sum\limits_{i=1}^{q_1+\dots+q_n}v_{i}
		\int_0^1
		\partial_i	F(\gamma v_{1}^u,\gamma v_2^u,\dots,\gamma v_{q_1+\dots+q_{n}}^u)\,   \,d\gamma,
	\end{array}
	\end{equation} 
		where $\partial_i$ here and in what follows denotes the partial derivative with respect to the $i$-th
	argument. Define
	$$
	\begin{array}{ll}
\displaystyle	I=\Biggl\{(s,r)\in\N^2\,:\, 1\le s\le n,\ 1+\sum\limits_{j=1}^{s-1}q_j\le r\le 
\sum\limits_{j=1}^{s}q_j,\\
	\qquad\qquad \qquad \displaystyle \int_0^1\partial_{r}{F}(\gamma v_{1}^u,\gamma v_2^u,\dots,\gamma v_{q_1+\dots+q_{n}}^u)\,d\gamma\ne 0
	\Biggr\},
	\end{array} 
	$$
	where the sum over the empty set equals zero. Note that the set $I$ is not empty, for else the representation 
	(\ref{QQ})--(\ref{baru}) is impossible and we immediately get a contradiction to 
	our assumption. Then,
	for an  arbitrarily fixed   $(s_0,r_0)\in I$,
	one can choose the initial function $\varphi$ such that
	$\varphi_{s_0}\left(x_{s_0r_0}\right)\ne 0$
	while $\varphi_s\left(x_{sr}\right)=0$
	for all other  $(s,r)\in I$. 
	On account of (\ref{baru}), the equality (\ref{Qphi1}) now reads
	$$
	u_j(x,t)=\varphi_{s_0}(x_{s_0r_0})\int_0^1\partial_{r_0} {F}(\gamma v_{1}^u,\gamma v_2^u,\dots,\gamma v_{q_1+\dots+q_{n}}^u)\,d
	\gamma\ne 0,
	$$
	contradicting 
		the FTS property of our problem. We, therefore,
	proved that the condition
	(\ref{q1}) is true for all $ t\in[T,2T]$.

Now we show that 	(\ref{q1}) is true for all $ t\ge 2T$.
To this end, observe that in the
autonomous case  the following formulas are true:
\begin{equation}\label{*1}
	\begin{array}{rcl}
	\omega_j(\xi,x,t+T)&=&\omega_j(\xi,x,t)+T,\quad t\ge 0,\\ [2mm]
	\left[Sv\right]_j(x,t)&=&c_j(x_j,x,t)v_j(\omega_j(x_j,x,T)+t-T),
	\quad t\ge T,
\end{array}
\end{equation}
for all  $v\in C(\R_+)^n$. Given $w\in C_h(\Pi)^n$, set 
 $z(x,t)=w(x,t+T)$. It follows that 
\begin{equation}\label{q2}
		\left[(SR)^q z\right](x,t)=\left[(SR)^q w\right](x,t+T),
	\quad t\ge T.	
\end{equation} 
Using the above argument for 
(\ref{q1}) for $t\in[T,2T]$
once again, we see that  $T>T_{opt}>1/a$.
On account of (\ref{*1}), we then have
	 $\omega_j(x_j(x,t),x,t)=\omega_j(x_j(x,t),x,t-T)+T>T$ for all $t>2T$,
$x\in[0,1]$, and $j\le n$.
Combining this with the FTS property,
we conclude that
$u(\cdot,t)=[Qu](\cdot,t)=[SRu](\cdot,t)\equiv 0$ for all  $t>2T$. 
Summarizing,
the condition
(\ref{q1}) stays true for all $t\ge T$, as desired.

Let $q$ be now chosen such that   (\ref{q1}) holds for $t\ge T$ and,
additionally,
the equality (\ref{tQ0}) is fulfilled in $\Pi^{3T}$.
Let  $w\in C_h(\Pi)^n$ be arbitrarily fixed.
Similarly to the proof of Theorem~\ref{control_inverse}, the function $[Q^{q}w](x,t)$ is a continuous solution
to 	(\ref{ss2d}), (\ref{ss11}), (\ref{ss5}),  (\ref{eq:h})  with $\varphi(x)=w(x,0)$
 in the domain
$\Pi^{3T}$. By (\ref{q1}), we  have
$\left[Q^{q}w\right](\cdot,T)\equiv 0$ and, hence
the function  $z^1(x,t)=\left[Q^{q}w\right](x,t+T)=\left[(SR)^{q}w\right](x,t+T)$ belongs to $ C_h(\Pi)^n$ and
is a continuous solution
	(\ref{ss2d}), (\ref{ss11}), (\ref{ss5}),  (\ref{eq:h})  with $\varphi(x)=0$
 in $\Pi^{2T}$. It follows  from (\ref{q1}) that
$$
0=\left[Q^{q}z^1\right](x,t)=
\left[(SR)^{q}z^1\right](x,t)\quad \mbox{ for } t\in[T,2T].
$$
Similarly to (\ref{q2}), we have
$$
	\left[(SR)^{q}z^1\right](x,t)=
\left[(SR)^{q}Q^qw\right](x,t+T)=\left[Q^{2q}w\right](x,t+T).
$$
Therefore, $\left[Q^{2q}w\right](\cdot,t)\equiv 0$ for $t\in[2T,3T]$.
In the next step we set $z^2(x,t)=\left[Q^{2q}w\right](x,t+2T)$. Due to the previous step, $z^2(\cdot,0)\equiv 0$ and, therefore, $z^2$
belongs to $C_h\left(\Pi\right)^n$ and is a continuous solution
to 	(\ref{ss2d}), (\ref{ss11}), (\ref{ss5}),  (\ref{eq:h})  with $\varphi(x)=0$
in $\Pi^{2T}$.  Similarly, for $t\in[T,2T]$, it holds $$0=\left[Q^{q}z^2\right](x,t)=\left[(SR)^{q}z^2\right](x,t)
=\left[(SR)^{q}Q^{2q}w\right](x,t+2T)=\left[Q^{3q}w\right](x,t+2T)
$$ 
 and, hence $\left[Q^{3q}w\right](\cdot,t)\equiv 0$ for $t\in[3T,4T]$. Proceeding further 
by induction, where on the $k$-th step we set $z^k(x,t)=\left[Q^{kq}w\right](x,t+kT)$, $k\ge 3$,
 we conclude that the desired condition (\ref{C00}) is true.
The proof of Case 1 is therewith complete. 

{\bf Case 2:} the problem (\ref{ss2d}), (\ref{ss11}), (\ref{ss5}),  (\ref{eq:h}) is FTS and all $L^2$-generalized solutions stabilize to zero in a finite time. Let $q$ be  as in Case 1. Using the same argument as in 
the proof of the necessity part of Theorem \ref{tm_1} in the same $L^2$-case, 
fix an arbitrary $w\in C_h(\Pi)^n$, put $\varphi(x)=w(x,0)$, and conclude that the function $
u=Q^qw\in C_h(\Pi)^n$ is a continuous solution    to
the problem (\ref{ss2d}), (\ref{ss11}), (\ref{ss5}),  (\ref{eq:h})
in the domain $ \Pi^{3T}$. 
Since any continuous solution is an $L^2$-generalized solution, 
then using the FTS property for the $L^2$-generalized solutions
and  (\ref{q1}),
we conclude that
$
\left[Q^qw\right](x,T)=0$
for all  $x\in[0,1]$. 
The proof is completed  by repeating the argument used at the end of  
Case~1.

\section{Examples}\label{ex}
\renewcommand{\theequation}{{\thesection}.\arabic{equation}}
\setcounter{equation}{0}
\subsection{Solving inverse problems}\label{inverse}
Let the boundary conditions (\ref{ss5}) be linear, namely 
\begin{equation}\label{n15l}
	u^{out}(t)=Pu^{in}(t),\quad t\ge 0,
\end{equation}
where $P=(p_{jk})$ is an  $n\times n$-matrix with constant entries.
We assume that  the matrix $P_{abs}=(|p_{jk}|)$ is nilpotent.
Then, due to \cite[Theorem~1.10]{KL2},
the  problem 
(\ref{ss2d}), (\ref{ss11}), (\ref{n15l}) is robust FTS, with
respect to   perturbations of the coefficients
$a_j$ and $b_{j}$.

Fix an arbitrary $r>0$ and consider the following abstract setting of the autonomous problem (\ref{ss2d}), (\ref{ss11}), (\ref{n15l}) on $L^2(0,1)^n$  (as studied, e.g., in \cite{Pril}, \cite{tikh18}):
\begin{eqnarray}
	\frac{d}{dt}u(t)={\A}u(t)+f ,\quad (0\le t \le r)
\label{n1}\\
u(0)=u_0,\quad
u(r)=u_{r},\label{fin}
 \end{eqnarray}
where  the operator
$ \A: D(\A)\subset L^2(0,1)^n\to L^2(0,1)^n$ is defined by 
$$
\begin{array}{rcl}
	\left(\A v\right)(x)&=&-A(x)v^{\prime} - B(x)v,\\ [1mm] D(\A)&=&\left\{v\in L^2(0,1)^n\,:\,v^{\prime}\in L^2(0,1)^n,\,v^{out}=
Pv^{in}\right\},
\end{array}
$$
and $u_0, u_{r}\in D(\A)$ are known functions.
Here $v^{out}$, $v^{in}$ are defined similarly to (\ref{v1}).
Solving the inverse problem  
(\ref{n1})--(\ref{fin}), we are looking for a
 couple of functions~$(u,f)$ such that $u\in C^1([0,r],L^2(0,1))^n$, $u(t)\in D(\A)$ for all $t\in [0,r]$, and $f\in L^2(0,1)^n$. 
 
Since the problem (\ref{n1})--(\ref{fin}) is autonomous, then, due to \cite[Theorem 2.3]{KmLyul},  the operator $ \A$ generates a $C_0$-semigroup~$S(t)$. Since the problem (\ref{n1})--(\ref{fin})
 is FTS, the semigroup $S(t)$ is nilpotent. Hence, there exists $T>0$ such that
 	$S(t)=0$ for all $t\ge T$. Accordingly to
 	\cite[Theorem 4]{tikh18}, for any
 $u_0, u_{r} \in D(\A)$, there is a unique function $f\in  L^2(0,1)^n$
 solving the inverse problem (\ref{n1})--(\ref{fin}). Moreover,
 this function admits the representation 
$$
f=\left\{
\begin{array}{rl}
	-{A}u_{r}  &\quad\mbox{if }\ r\ge T\\ 
\displaystyle	-{A}u_{r}+A\sum_{k=1}^{n_0} S(kr)(u_0-u_{r}) &\quad\mbox{if }\ r<T,
\end{array}
\right.
$$
where 
 $n_0=\lceil T/r\rceil-1$. Recall that  $\lceil x\rceil$ denotes the integer nearest to  $x$ from above. The unknown function
  	$u(t)$ is then given by the formula
$$
u(t)= S(t)u_0+\int_0^t S(s)f\,ds,\quad 0\le t\le r.
$$   

Now, using Theorem \ref{control_inverse},
we conclude that  there exists $k=k(T)\in\N$ such that for all $x\in[0,1]$
it holds that 
$$
[S(t)u_0](x)=\left\{
\begin{array}{rl}
	\left[Q^k w\right](x,t)  &\quad\mbox{if }\  t\le T\\ 
	0 &\quad\mbox{if }\ t> T,
\end{array}
\right.
$$
the formula being true for any $w\in C_h(\Pi)^n$ such that
   $ w(x,0)=u_0(x)$.

\subsection{Nonlinear boundary conditions and  FTS property} In the domain $\Pi$ we consider the $2\times 2$-decoupled system 
	\begin{equation}\label{1li111}
		\partial_tu_1+\partial_xu_1=0,\quad \partial_tu_2-\partial_xu_2=0
	\end{equation}
	with the nonlinear boundary conditions
	\begin{equation}\label{1li211}
		u_1(0,t)=r(t) \sin(u_2(0,t)),\quad u_2(1,t)=\sin^2 (s(t) u_1(1,t))\end{equation}
	and the initial conditions 
	\begin{equation}\label{eq:in1}
		u_1(x,0)=\varphi_1(x),\quad u_2(x,0)=\varphi_2(x).
	\end{equation}
	Here $r$ and $s$ are smooth and uniformly bounded functions for
	$t\ge 0$. Note that the  boundary conditions are of the type    (\ref{eq:h1}).
	Our aim is, using Theorem \ref{tm_1}, to 
	find conditions on the functions  $r$ and $s$ such that the problem
	(\ref{1li111})--(\ref{eq:in1}) is  FTS. 
	
	The operator $Q$ defined by 
	(\ref{Q}) is now specified to
	\begin{eqnarray*}
		[Qu]_1(x,t)&=&\left\{
		\begin{array}{ll}
			\varphi_1(x-t)&\mbox{if}\ x> t \\
			r(t-x) \sin(u_2(0,t-x))&\mbox{if}\ t-x\ge 0,
		\end{array}
		\right.\\ [2mm]
		[Qu]_2(x,t)&=&\left\{
		\begin{array}{ll}
			\varphi_2(x+t)&\mbox{if}\ t+x<1 \\
			\sin^2 (s(t+x-1) u_1(1,t+x-1)) &\mbox{if}\ t+x\ge 1.
		\end{array}
		\right.
	\end{eqnarray*}
	The second power of $Q$ is then given by
	\begin{eqnarray*}
		[Q^2u]_1(x,t)&=&\left\{
		\begin{array}{ll}
			\varphi_1(x-t)&\mbox{if}\ x> t \\
			r(t-x)\sin(\varphi_2(t-x))&\mbox{if}\ 0\le t-x<1\\
			r(t-x)\sin\left(\sin^2 (s(t-x-1)u_1(1,t-x-1))\right)&\mbox{if}\ 1\le t-x,
		\end{array}
		\right.\\
		\left[Q^2u\right]_2(x,t)&=&\left\{
		\begin{array}{ll}
			\varphi_2(x+t)&\mbox{if}\ t+x<1 \\
			\sin^2 \left(s(t+x-1) \varphi_1(2-(t+x))\right)&\mbox{if}\ 1\le t+x<2\\
			\sin^2 \left(s(t+x-1)r(t+x-2)\sin(u_2(0,t+x-2))
			\right) &\mbox{if}\ 2\le t+x.
		\end{array}
		\right.\\
	\end{eqnarray*}
			It follows that if there exist reals  $T_1>0$ and $T_2>0$
		with
		\begin{equation}\label{suf2}
			T_2-T_1\ge 1\quad\mbox{and}\quad	\Bigl(r(t)=0 \, \mbox{ and }\ s(t)=0\,  \mbox{ for }\ T_1\le t\le T_2\Bigr),
		\end{equation}
		then  the condition (\ref{C0}) is true  with $k=1$. If there exist reals  $T_1>0$ and $T_2>0$ with
		\begin{equation}\label{suf1}
			T_2-T_1\ge 2\quad\mbox{and}\quad\Bigl(r(t)=0 \, \mbox{ or } s(t)=0 \, \mbox{ for } \, T_1\le t\le T_2\Bigr),
		\end{equation}
		then  the condition (\ref{C0}) is true  with $k=2$.			
		In other words, (\ref{suf2}) and (\ref{suf1}) are two sufficient conditions for the problem 
		(\ref{1li111})--(\ref{eq:in1}) to be FTS.

					\subsection{Theorem  \ref{tm_1} does not extend  for nonhomogeneous  boundary conditions}  In the domain $\Pi$, we consider  the  $2\times 2$-decoupled system
			(\ref{1li111}) with the initial conditions (\ref{eq:in1}) and
			the boundary conditions
			\begin{equation}\label{lik}
				u_1(0,t)=g(t),\quad u_2(1,t)=u_1(1,t).
			\end{equation}
			
			Fix    $g$ to be a smooth bounded function such that
			\begin{eqnarray*}
				g(t)=\left\{
				\begin{array}{llc}
					0&\quad\mbox{if }\ 0\le t\le 4 \\
					\neq 0\ &\quad\mbox{if }\ 4<t.
				\end{array}
				\right.
			\end{eqnarray*}
			The formula (\ref{Q})  then reads 
			\begin{eqnarray*}
				\begin{array}{ll}
					[Qu]_1(x,t)=\left\{
					\begin{array}{rl}
						\varphi_1 (x-t)&\mbox{if}\ x> t \\
						g(t-x) &\mbox{if}\ t-x\ge 0,
					\end{array}
					\right.\,\,
					[Qu]_2(x,t)=\left\{
					\begin{array}{rl}
						\varphi_2(x+t)&\mbox{if}\ t+x<1 \\
						u_1(1,t+x-1) &\mbox{if}\ t+x\ge 1,
					\end{array}
					\right.
				\end{array}
			\end{eqnarray*}
			implying that
			\begin{eqnarray*}
				\left[Q^2u\right]_1(x,t)=\left[Qu\right]_1(x,t),\quad
				\left[Q^2u\right]_2(x,t)&=&\left\{
					\begin{array}{ll}
						\varphi_2(x+t)&\mbox{if}\ t+x<1 \\
						\varphi_1(2-(t+x)) &\mbox{if}\ 1\le t+x< 2\\
						g(t+x-2) &\mbox{if}\ 2\le t+x.
					\end{array}
					\right.
				\end{eqnarray*} 
				It follows that $
				\left[Q^2u\right](x,3)\equiv 0,
				$
				while the problem  (\ref{1li111}), (\ref{eq:in1}), (\ref{lik}) is not FTS.

\section*{Acknowledgement}
This work was supported by the VolkswagenStiftung Project ``From Modeling and Analysis to Approximation''.


\begin{thebibliography}{99.}%

\bibitem{AbdCos02}
Amato, F., Ariola  M., Abdallah C., Cosentino, C.: Application  of 
finite-time stability concepts to the control of ATM networks.   Proceedings: Allerton Conference on Communication, Control, and Computing {\bf 40(2)}, 1071--1079 (2002).   



\bibitem{CosTom14} Amato, F.,  Ambrosiano, R., Ariola, M., Cosentino, C., Tommasi, G.D., et al.: Finite-time stability and control. Springer, Berlin (2014).

\bibitem{Or105} Orlov, Y.:
Finite Time Stability and Robust Control Synthesis of Uncertain Switched Systems.  SIAM J.  Control Optim. {\bf 43(4)}, 1253--1271 (2004).




\bibitem{corbas} Bastin, G.,  Coron, J.-M.:  Stability and boundary stabilization of 1-d hyperbolic systems. In:
 Progr. Nonlinear Differential Equations and Their Appl., \textbf{  88}, Birkh\"aser (2016).
 
 


\bibitem{gugat} Gugat, M.: Optimal boundary control and boundary stabilization of hyperbolic systems. Birkh\"aser, Basel (2015).


\bibitem{hege} He, W., Ge, S.-Z.:  Robust adaptive boundary control of a vibrating string under unknown time-varying disturbance.  IEEE Transactions on Control Systems Technology. \textbf{20(1)}, 48--58 (2012).


  \bibitem{hegez}He, W., Ge, S.-Z., Zhang, S.: Adaptive boundary control of a flexible marine installation system. Automatica. \textbf{47(12)}, 2728--2734 (2011).


\bibitem{BerKrs} Karafyllis, I., Krstic, M.: Input-to-state stability for PDEs. Springer,  (2019).


\bibitem{ijdsde}
Kmit, I.: Classical solvability of nonlinear initial-boundary problems for
first-order hyperbolic systems. Intern. J. Dynamic Systems Different.
	Equat. \textbf{1(3)}, 191--195 (2008).

\bibitem{KL2} Kmit, I., Lyul`ko, N.: Finite time stabilization of nonautonomous first-order hyperbolic systems. SIAM J. Control Optim. \textbf{59(5)}, 3179--3202 (2021).

\bibitem{lyul_I}Lyul`ko, N.: 1D hyperbolic systems with  nonlinear boundary conditions, I:
$L^2$-generalized solutions. To appear in {\it Research Perspectives:
Analysis, Applications, and Computations - Selected contributions of the 13th ISAAC Congress, Ghent, Belgium, 2021}, Basel: Birkh\"auser.

\bibitem{KmLyul} Kmit, I., Lyul`ko, N.: Perturbations of superstable linear hyperbolic systems.
 J. Math. Anal. Appl. \textbf{460(2)}, 838--862 (2018).

\bibitem{lakra}
Pavel, L.:
Classical solutions in Sobolev spaces for a class of hyperbolic Lotka-Volterra systems. SIAM J. Control Optim. \textbf{51(3)}, 2132--2151 (2013).



\bibitem{per}Perruquetti, A., Barbot, J.P.: Sliding Mode Control in Engineering.  	New York: M. Dekker
 (2002).

\bibitem{Pril} 
Prilepko,  A.I., Orlovsky D.G., Vasin I.A.: Methods for solving inverse problems in matthematical physics. Taylor \& Francis (2000).



\bibitem{tikh18} Tikhonov I., Vu Nguyen Son Tung:
The solvability of the  inverse problem for the evolution equation with a superstable semigroup.  RUDN Journal of MIPh.  \textbf{26(2)}, 103--118 (2018).


 
\bibitem{udw05} Udwadia, F.E.: Boundary control, quiet boundaries, super-stability and super-instability.
 Appl. Math. Comput. \textbf{164(2)}, 327--349 (2005).

\bibitem{udw12} Udwadia, F.E.: On the longitudinal  vibrations of a bar with  viscous  boundaries: super-stability, super-instability
and loss damping.
 Int. J. Eng. Sci. \textbf{50(1)}, 79--100 (2012).


 \bibitem{xu} Xu, G. Q.: Stabilization of string system with linear boundary
feedback. Nonlinear Analysis: Hybrid Systems.  \textbf{1}, 383--397 (2007).



 
 


\end{thebibliography}
\end{document}